\documentclass[12pt,twoside]{article} 
\usepackage{amsmath, amsfonts, amsthm, enumerate}
\newcommand {\C}{\mathbb C}
\newcommand{\R}{\mathbb R}
\newcommand{\B}{\mathcal B}
\newcommand{\Oc}{\mathcal O}

 \newtheorem {thm}{Theorem}
 \numberwithin{thm}{section}
 \newtheorem {lem}[thm] {Lemma}
 \newtheorem{prop}[thm]{Proposition}
 \numberwithin{equation}{section}
\begin{document}
 
 \begin{titlepage}
 \title{ \bf An Embedding Theorem for Pseudoconvex Domains in Banach Spaces\thanks{This research supported in part by NSF grant DMS-0203072 and a Bilsland Dissertation Fellowship from the Graduate School of Purdue University.  This paper is the result of my thesis research under the direction of L\'aszl\'o Lempert, and I would like to thank him for his patient and wise guidance.  I would also like to thank Andreas Defant and David P\'erez-Garc\'{\i}a for their suggestions regarding tensor products and directing me to the paper [GR], and Steve Bell for his suggestions. }}
 \author{Aaron Zerhusen\\ Department of Mathematics\\ Purdue University\\ West Lafayette, IN 47907-2067, USA\\ azerhus@math.purdue.edu}
 
 \end{titlepage} 
\date{}
 \maketitle
 \abstract
 We show that a pseudoconvex open subset of a Banach space with an unconditional basis is biholomorphic to a closed direct submanifold of a Banach space with an unconditional basis.  \endabstract
 \section{Introduction}
A fundamental result in complex geometry, shown by Bishop, Narasim\-han and Remmert, is that an $n$-dimensional Stein manifold is biholomorphic to a closed submanifold of $\C^{2n+1}$ [B, N, Re].  One reason that this is of such importance is in the study of holomorphic functions on manifolds.  Cartan's Theorem B implies that a holomorphic function on a closed submanifold $M\subset \C^k$ extends to a function holomorphic on all of $\C^k$.  Therefore, once a given Stein manifold is embedded in some $\C^k$, its holomorphic functions can be studied through holomorphic functions on $\C^k$.   In [L2] Lempert showed a result similar to the Bishop-Narasimham-Remmert embedding theorem for Banach spaces, namely, if $\Omega $ is a pseudoconvex domain in a Banach space $X$, and $X$ has an unconditional basis, then $\Omega $ is biholomorphic to a closed direct submanifold of a separable Banach space $Y$.  However, his proof does not guarantee that $Y$ would have an unconditional basis or Hilbert space structure if $X$ does.  This is an important issue, since many complex analytic results, such as approximation theorems, rely upon the existence of an unconditional basis.  In particular, in [L3] and [P] it is shown that if $M$ is a closed direct submanifold in $Y$, a Banach space with unconditional basis, then holomorphic functions on $M$ extend to holomorphic functions on $Y$.  In this paper we resolve the above issue by improving upon Lempert's embedding theorem. 
Our main theorem is 
\begin{thm}
Let $X$ be a Banach space with an unconditional basis, and $\Omega \subset X$ pseudoconvex and open.  Then $\Omega$ is biholomorphic to a closed direct submanifold of a Banach space $Y$, which has an unconditional basis.  
Further, if $1<p<\infty$ and $X$ is a separable $L^p$ space, then $Y$ can also be taken to be a separable $L^p$ space.  
\end{thm}
Just as in [L2], this is an immediate consequence of a domination theorem.  

\begin{thm} 
Let $X$ be a Banach space with an unconditional basis, $\Omega \subset X$ pseudoconvex and open, and $u:\Omega \to \mathbb R$ locally bounded.  There is a Banach space $V$ with unconditional basis and a holomorphic $f: \Omega \to V$ such that $u(x)\leq \|f(x)\|$ for $x\in \Omega$.  If $1<p<\infty$ and $X$ is a separable $L^p$ space, then $V$ may also be taken to be a separable $L^p$ space.  
\end{thm}

\begin{proof}[Proof of Theorem 1.1 from Theorem 1.2] Let $d:\Omega \to \mathbb R$, $d(x)= \text{dist}(x, \partial \Omega) $, and let $f\in \mathcal O(\Omega , V) $ such that $\|f(x)\| \geq \frac{1}{d(x)} $.  Then the graph of $f$ is a closed direct submanifold of $X\times V$ which is biholomorphic to $\Omega$.  This proves the theorem, since if $X$ and $V$ have unconditional bases so does $X\times V$, and if $X$ and $V$ are $L^p$ spaces so is $X\times V$. \end{proof}

\section{Definitions}
We will require some notions from the theory of Banach spaces.  A Schauder basis for a Banach space $X$ (for this paper, we always mean over $\C$) is a countable sequence of vectors $e_1, e_2, \ldots \in X$ such that  for every $x\in X$ there is a unique sequence of complex numbers $\lambda _1, \lambda_2, \ldots \in \C$ such that 
\[x=\sum \lambda _ne_n\] where the sum converges in norm.  Having a Schauder basis implies that the space is separable.  If the sequence converges in norm after arbitrary rearrangements, it is called an unconditional basis, or equivalently, if $\sum_1^\infty \lambda_n \theta_n e_n$ converges whenever $\sum_1^\infty \lambda_n  e_n$ does, for all choices of $\theta _n =\pm1$.  Following [LT], for $x=\sum_1^\infty \lambda_n e_n$ and every $\theta=\{\theta _n\}_{n=1}^\infty, \theta_n=\pm 1$ define a linear operator $M_\theta x= \sum_1^\infty \lambda_n \theta_n e_n$.  The principle of uniform boundedness guarantees that $c=\sup_\theta \|M_\theta\|<\infty$, and we call $c$ the unconditional constant of the basis.  We call a basis with unconditional constant 1 a 1-unconditional basis.  Upon renorming with an equivalent norm, any space with an unconditional basis has a 1-unconditional basis [LT, p. 18].  The usual basis of $l^p$ for $1\leq p<\infty$ and the Haar basis of $L^p[0,1]$ for $1<p<\infty $ are examples of unconditional bases  (although the Haar basis is not 1-unconditional).  $L^1[0,1]$, however does not have an unconditional basis.  Since any separable $L^p$ space is isometrically isomorphic to $l^p_n,\ l^p,\ L^p[0,1]$, or a direct sum of these, for $1<p<\infty$ any separable $L^p$ space has an unconditional basis.  With a choice of basis $\{e_n\} $ for a Banach space $X$ we associate a set of projections.  Let $\pi_N:X\to X$ be the projection $\pi_N \sum_{n=1}^\infty \lambda_n e_n=\sum_{n=1}^N\lambda_n e_n$ and $\rho_N=$id$-\pi_N$.  If $ \|\pi_N\|=1$ for all $N$, we say that $\{e_n\}$ is monotone.  It is immediate that every 1-unconditional basis is monotone, and it is also true that the Haar basis is a monotone basis for $L^p(0,1)$ for all $1\leq p<\infty$ [LT, p. 3].  Therefore, for all $1< p <\infty$ any separable $L^p$ space has a monotone unconditional basis.

We will also require certain notions from the theory of complex analysis in Banach spaces.  For more on this subject, see [D or M].  For an open set $\Omega\subset X$ to be pseudoconvex means that $\Omega \cap Y$ is pseudoconvex in the usual sense for all finite dimensional subspaces $Y\subset X$.  A Hausdorff space $N$ is a complex manifold modeled on a Banach space $Z$ if $N$ has an open cover $\{U_\alpha\}$, for each $\alpha$ we specify a homeomorphism $\phi _\alpha$ from $U_\alpha$ to an open subset of $Z$, and $\phi_\alpha \phi^{-1}_\beta$ are holomorphic where defined.  For $M\subset N$ to be a complex submanifold means that as a pair $(M,N)$ is locally biholomorphic to $(X,Z)$, $X$ a Banach subspace of $Z$.   If $X$ is a complemented subspace of $Z$ we say that $M$ is a direct submanifold of $N$.  

\section{The Proof of Theorem 1.2}

Let $\B$ denote the covering of $\Omega \subset X $
\[ \B=\B_\Omega =\{ \text{balls }B : \bar{B}\subset \Omega,\  2 \text{diam }B < \text{diam }\Omega \}.\]
\begin{prop}
Suppose either that (a) X is a Banach space with a 1-unconditional basis, or (b) $1<p<\infty$ and $X$ is a separable $L^p$ space.  Given $\Omega \subset X$ pseudoconvex open, and $u:\Omega \to \R$, suppose that for all $B \in \B_\Omega$ there are a Banach space  $V_B$ with 1-unconditional basis in case (a), separable $L^p$ space in case (b), and $f_B \in \Oc (B; V_B)$ such that $u(x)\leq \|f_B(x)\|_{V_B}$ for all $x\in B$.  Then there are a Banach space $V$ and $f\in \Oc (\Omega, V)$ such that in case (a) $V$ has 1-unconditional basis, in case (b) $V$ is a separable $L^p$ space, and $u(x)\leq \|f(x)\|_V$ for all $x\in\Omega$. 
\label{vb}\end{prop}

\begin{proof}[Proof of Theorem 1.2 from Proposition 3.1.]  This proof is the same as in [L2].  Let $u:\Omega \to \R$ be locally bounded and assume that on $\Omega$ holomorphic domination of $u$ with values in a space with 1-unconditional basis (respectively, a separable $L^p$ space) is not possible.  Then, by Proposition 3.1,  there exists a ball $B_1\subset \Omega$ such that diam $B_1 < 1/2$ diam $\Omega$, on which holomorphic domination  is not possible.  By again applying Proposition 3.1, but with $B_1$ instead of $\Omega$, there exists a ball $B_2$, $\bar{B_2} \subset B_1$ with diam $B_2<$ diam $B_1/2$, on which $u$ cannot be dominated by a holomorphic function.  Continuing  in this way we get a descending sequence of balls on which holomorphic domination is impossible.  Let $x=\bigcap_1^\infty B_j$.  Thus holomorphic domination is not possible on any neighborhood of $x$, contradicting the local boundedness of $u$. 
\end{proof}

\section{Exhaustion techniques}
The rest of this paper is devoted to proving Proposition 3.1.  Let $X$ be as there; $e_1, e_2, \dots $ denote, in case (a), a 1-unconditional basis; in case (b), a monotone unconditional basis.  Recall the projections $\pi_N$ and $\rho_N$ from Section 2.

To prove Theorem 1.2 or Proposition \ref{vb} in the finite dimensional setting, a reasonable technique would be to exhaust $\Omega$ by compact sets.  In the setting of an infinite dimensional Banach space, this is not possible as any compact set has empty interior.  Lempert proposes instead exhausting by ball bundles of a certain type.  The definitions and theorems in this section are from [L2].
In order to exhaust $\Omega$, we consider a special type of ball bundle.  Let 
\begin{equation} d(x)=\min \{1, \text{dist}(x,X\setminus \Omega)\}, \label{d} \end{equation}
 and $0<\alpha<1$.  For any positive integer $N$, define
\[ D_N\langle \alpha \rangle = \{t\in \pi_NX : \|t\|<\alpha N, 1<\alpha Nd(t)\}\]
\[ \Omega _N\langle\alpha \rangle = \{x\in X: \pi_Nx \in D_N\langle \alpha \rangle, \|\rho_Nx\|<\alpha d(\pi_Nx)\}. \]

\begin{thm}
\begin{enumerate}[a.]
  \item Each $\Omega_N\langle \alpha \rangle$ is pseudoconvex.  
  \item $\overline{ \Omega _n \langle \gamma \rangle } \subset \Omega _N \langle \beta \rangle $ if $n\leq N, \gamma \leq \beta/4$.
 
  \item For fixed $\gamma$ each $x\in\Omega$ has a neighborhood that is contained in all but finitely many $\Omega_N\langle \gamma \rangle$.
\end{enumerate}
\end{thm}
This is proven in [L2, Proposition 3.1].  Note that in the proof of (b) only the monotonicity of the basis is used.

Just as in many situations in the theory of finitely many variables, we need a result regarding holomorphic approximation on certain sets.
\begin{lem} \label{approx}  Assume $X$ has an unconditional basis with basis constant $c$.  If $\gamma<2^{-7}\alpha /c$ and $V$ is a Banach space, then any $\psi \in \Oc (\Omega_N\langle \alpha \rangle ;V)$ can be approximated by $\phi \in \Oc (\Omega ;V)$, uniformly on $ \Omega_N\langle \gamma \rangle$.
\end{lem}
We first consider the simpler case of balls about the origin.  Let $B(\mu)=\{x\in X:\|x\|<\mu\}$.  
\begin{lem}Assume $X$ has an unconditional basis with basis constant $c$.  Then for any Banach space $(W,\|\ \|_W)$, $\epsilon>0$, and $g\in\Oc(B(1);W)$ there is an $h\in\Oc(X;W)$ that satisfies $\|g-h\|<\epsilon$ on $B(1/(2c))$.  
\end{lem}
\begin{proof}
It was shown by Lempert in [L1] and by Josefson in [J] that every space with an unconditional basis allows such approximation on a ball when $c=1$.  Let $\|x\|'=\sup_\theta \|M_\theta x\|$.  This is a norm equivalent to $\|\ \|$ under which $e_1, e_2, \dots$ is a 1-unconditional basis.  Let $B'(\mu)=\{x\in X:\|x\|'<\mu\}$.  For all $x\in X$ we have the inequality 
\[ \|x\|\leq \|x\|'\leq c\|x\|,\]
which gives that $B(1/(2c))\subset B'(1/2)$ and $B'(1)\subset B(1)$.  The lemma is therefore an immediate consequence of the Josefson-Lempert approximation results.  
\end{proof}

Now the proof of  Lemma \ref{approx} is a special case of [L2, Theorem 3.3] with $\mu=1/(2c)$.

\section{Tensor products and sums of Banach spaces}

In the proof of Proposition 3.1, we will use tensor products of Banach spaces.  For the general theory of tensor products see [DF or Ry].   If $X$ and $Y$ are Banach spaces, as usual $X\otimes Y$ will denote the tensor product of the underlying complex vector spaces.  We will use $X \overline{\otimes}  Y$ to denote the completion of $X\otimes Y$ with respect to a certain norm which preserves key features of $X$ and $Y$ in the following sense.  If $X=L^p(\mu)$ and $Y=L^p(\nu)$, we define the norm of $\sum_1^n x_j\otimes y_j \in X\otimes Y$ by 
\begin{equation} 
	\| \sum_1^n x_j\otimes y_j \|=\bigl(\int |\sum_1^n x_j(s) y_j(t)|^p d\mu(s) d\nu(t) \bigr)^{1/p}.
	\end{equation}
In particular, Fubini's theorem gives
\begin{equation}
	 \|x\otimes y\|=\|x\|\|y\|. 
 \label{crossnorm}	 \end{equation}
Under this norm $X\overline{\otimes}Y$ is isometrically isomorphic to $L^p(\mu\times\nu)$ [DF p.79].

If $X$ and $Y$ are spaces with 1-unconditional bases $(e_n)$ and $(f_n)$ respectively, Grecu and Ryan define a norm on $X\otimes Y$ with the property that $(e_n\otimes f_m)$ is a 1-unconditional basis for the completion $X\overline{\otimes}Y$~[GR].  How they define this norm and other properties it has are not of immediate importance to us, just that it preserves unconditional structure and also satisfies (\ref{crossnorm}).

We now look at certain facts about sums of Banach spaces.  If $(W_m)$ is a sequence of Banach spaces, each with 1-unconditional basis, or each an $L^p$ space, let
\[ W=  \overline{\bigoplus_{m=0} ^\infty} W_m\]
where the bar represents the completion with respect to the $l^p$ norm if $W_m$ are $L^p$ spaces or the $l^1$ norm if they are spaces with 1-unconditional basis.  Concretely, 
\[W=\{w=(w_m):w_m\in W_m, \|w\|=(\sum_{m=0}^\infty \|w_m\|^p)^{1/p}<\infty.\}\]  Clearly, $W$ is a separable Banach space.
\begin{lem}
	\begin{enumerate}
	\item If $W_m$ all have 1-unconditional bases, so does $W$.
	\item If $W_m$ are all $L^p$ spaces, so is $W$.
	\end{enumerate}
\end{lem}
\begin{proof}	Let $(e_{mn})_{n=1}^\infty$ be a 1-unconditional basis for $W_m$, $\theta =(\theta_{mn}), \theta_{mn}=\pm 1$.  If $x\in W$, $x=\sum \lambda_{mn} e_{mn}$.  
\begin{align*}
	\| \sum \theta _{mn} \lambda _{mn} e_{mn}\| &=(\sum_{m=0}^\infty \|\sum_{n=0}^\infty \theta_{mn} \lambda_{mn} e_{mn} \|) \\
	 &\leq (\sum_{m=0}^\infty \|\sum_{n=0}^\infty  \lambda_{mn} e_{mn} \|) \\
	 &\leq \|x\| \\
\end{align*}So $(e_{mn})$ is a 1-unconditional basis for $W$.  If $W_m$ are $L^p$ spaces, $W$ is isometrically isomorphic to the corresponding $L^p$ space on the disjoint union of the appropriate measure spaces.  
\end{proof}

\section{The proof of Proposition \ref{vb}}

Fix $N\in \mathbb N$, and let $\pi =\pi_{N+1}$ and $\rho=\rho_{N+1}$.   If $A\subset \pi X$ and $r:A\to \R$ is continuous, then ball bundles of the form 
\[A(r)=\{x\in X : \pi x\in A, \|\rho x\| <r(\pi x)\} \]
\[A[r]=\{x\in X : \pi x\in A, \|\rho x\| \leq r(\pi x)\} \]
are called sets of type B.  Note that the sets $\Omega_N\langle\alpha \rangle$ are specific examples of sets of type B.

\begin{lem}
Assume either that (a) $X$ is a Banach space with 1-uncondi\-tion\-al basis or (b) $1<p<\infty$ and X is a separable $L^p$ space.  Let  $\Omega\subset X$ be open and $A_2\subset \subset A_3\subset\subset A_4$ be relatively open subsets of $\Omega \cap \pi X$, $A_1\subset A_2$ compact and plurisubharmonically convex in $A_4$.  Let $r_i:A_4\to (0,\infty)$ be continuous, $1\leq i\leq 4$, and $r_i<r_{i+1}$, with $-\log r_1$ plurisubharmonic.  Assume that any Banach space valued holomorphic function on $A_4(r_4)$ can be uniformly approximated on $A_3(r_3)$ by functions holomorphic on $\Omega$.  Then there exist a Banach space  $W$ and a function $\psi\in \Oc(\Omega, W)$ such that in case (a) $W$ has a 1-unconditional basis, in case (b) $W$ is a separable $L^p$ space, and
\begin{enumerate}
	\item $\|\psi(x)\| \leq 1/2$ if $x\in A_1[r_1]$
	\item $\|\psi(x)\| \geq 2$ if $x\in A_3(r_3) \setminus A_2(r_2)$.
\end{enumerate}
\end{lem}

\begin{proof}
Define Hartogs sets in $\pi X \times \C \approx \C^{N+2}$:
\[H_1=\{ (s,\lambda)\in A_1\times \C : |\lambda|\leq r_1(s)\}\]
\[H_i=\{ (s,\lambda)\in A_i\times \C : |\lambda|< r_i(s)\} , 2\leq i\leq 4.\]
The set $H_1$ is plurisubharmonically convex, hence holomorphically convex, in $H_4$ [H].  Since $\overline H_3\setminus H_2$ is compact, there exist finitely many holomorphic functions $\phi_j \in \Oc (H_4)$, $j=1,\dots, J$, such that if $z\in H_1$ we have $|\phi_j(z)|<1/4$ for all $j$ and if $z\in H_3\setminus H_2 $ there is a $j$ such that $|\phi_j(z) |>4$.  There exists a positive $\eta<1$ such that with 
\[H'_1=\{(s, \lambda)\in A_1\times \C: |\lambda|\leq \eta^{-1} r_1(s)\}\]
\[H'_2=\{ (s,\lambda)\in A_2\times \C : |\lambda|< \eta r_i(s) \} , \]
$|\phi_j(z)|< 1/2$ for all $j$ if $z\in H'_1$ and $|\phi_j(z)|>2$ for some $j$ if $z\in H_3\setminus H'_2$.  By replacing $\phi_j$ with $\phi_j ^n$ for a suitably large $n$ we can assume that $|\phi _j (z)|<(1-\eta)/(4J)$ for all $j$ if $z\in H'_1$ and $|\phi _j (z)|>4/(1- \eta)$ for some $j$ for each  $z\in H_3\setminus H'_2$.   Let $W'=\overline\oplus _{n=0}^\infty X^{\overline \otimes n}$, taking $X^{\overline \otimes 0}$ to be $\C$.  As shown before, $W'$ has a 1-unconditional basis if $X$ does and $W'$ is a separable $L^p$ space if $X$ is.  Since $\phi _j$ is holomorphic, $\phi_j(s,\lambda)=\sum_{n=0}^\infty a_{nj}(s) \lambda^n$, where $a_{nj}$ depends holomorphically on $s\in A_4$ and locally on $H_4$ convergence is absolute and uniform.  Define $\psi _j$ from $A_4(r_4)$ to $W'$ as 
\begin{equation}
	\psi_j(x)=(a_{nj}(\pi x)(\rho x)^ {\otimes n} )_{n=0}^\infty.\end{equation}  By the absolute and uniform convergence of the power series of $\phi_j$, $\psi_j$ is holomorphic.  

For any $x\in A_4(r_4)$ and any $j=1,\dots, J$, by definition of the norm of $W'$ (and since the $l^1$ norm dominates the $l^p$ norm for $p\geq1$) 
\[ \|\psi_j(x)\| \leq \sum_{n=0}^\infty |a_{nj}(\pi x)| \|\rho x\|^n.\]
If $x\in A_1[r_1]$,
\[  \|\psi_j(x)\| \leq \sum_{n=0}^\infty |a_{nj}(\pi x)| (r_1(\pi x))^n.\]
Fixing $s\in A_1$ and using Cauchy estimates for $\psi _j$ on $|\lambda |=\eta^{-1} r_1(s)$ we see that $|a_{nj}(s)|\leq (1-\eta )( \eta /r_1(s))^{n}/(4J)$, so 
\[ \|\psi_j(x)\| \leq (1-\eta)/(4J) \sum_{n=0}^\infty  \eta ^{n}=1/(4J).\]

Next let $x\in A_3(r_3)\setminus A_2(r_2)$ so that $(\pi x, \|\rho x\|)\in H_3\setminus H_2$.  In case (a), $W'$ is an $l^1$ sum and it is immediate that $\|\psi_j(x)\|\geq |\phi_j(\pi x, \|\rho x\|)|$, hence $\|\psi_j(x)\|>4$ for some $j$.  In case (b), i.e. $X$ is an $L^p$ space, choose $j$ so that $|\phi_j(\pi x, \eta \|\rho x \|)|>4/(1-\eta)$.  By H\"older's inequality: 
\begin{align*} \|\psi_j(x)\| & =( \sum_{n=0}^\infty( |a_{nj}(\pi x)| \|\rho x\|^n)^p)^{1/p}  \\
	&\geq (\sum_{n=0}^\infty |a_{nj}(\pi x)| \| \eta \rho x\|^n) (\sum_{n=0}^\infty \eta ^{nq})^{-1/q}\\
 	&\geq |\phi_j(\pi x, \eta\| \rho x\|)| (\sum_{n=0}^\infty \eta ^{n})^{-1}\\
	&>4.\end{align*}
Define $W=\oplus_{j=1}^JW'$ and $\tilde{\psi}(x)=(\psi_1(x), \dots \psi_J(x))$.  Let $\psi\in \Oc(\Omega, W)$ approximating $\tilde {\psi}$ within $1/4$ on $A_3(r_3)$.  Then $ \| \psi(x)\| \leq 1/2$ if $x\in A_1[r_1]$, and $\|\psi(x)\| \geq 2$ if $x\in A_3(r_3) \setminus A_2(r_2)$.  

\end{proof}

In the next lemma the geometric assumptions are the same as in the previous one.

\begin{lem}
Assume either that (a) $X$ and $Y$ are Banach spaces with 1-unconditional basis or (b) $1<p<\infty$ and  $X$ and $Y$ are separable $L^p$ spaces.  Let $\Omega\subset X$ is open and $A_2\subset \subset A_3\subset\subset A_4$ be relatively open subsets of $\Omega \cap \pi X$, $A_1\subset A_2$ compact and plurisubharmonically convex in $A_4$.  Let $r_i:A_4\to (0,\infty)$ be continuous, $1\leq i\leq 4$ and $r_i<r_{i+1}$ with $-\log r_1$ plurisubharmonic.  Assume that any Banach space valued holomorphic function on $A_4(r_4)$ can be uniformly approximated on $A_3(r_3)$ by functions holomorphic on $\Omega$.  Let $\epsilon >0$.   Then given a function $g\in \Oc (X,Y)$, there exist a Banach space $Z$ and a function $h\in \Oc(\Omega, Z)$ such that in case (a) $Z$ has a 1-unconditional basis, in case (b) $Z$ is a separable $L^p$ space, and 
\begin{enumerate}
	\item $\|h(x)\| \leq \epsilon$ if $x\in A_1[r_1]$
	\item $\|h(x)\| \geq \|g(x)\|$ if $x\in A_3(r_3) \setminus A_2(r_2)$.
\end{enumerate}
\end{lem}

\begin{proof}
Without loss of generality, assume $\|g(x)\|\geq \epsilon$ for all $x\in A_4(r_4)$.
By Lemma 6.1, there are a Banach space $W$ with 1-unconditional basis (or an $L^p$ space, as appropriate) and a function $\psi \in \Oc (\Omega, W)$  with $\|\psi(x)\| \leq1/2$ for $x\in A_1[r_1]$ and $\|\psi(x)\|\geq 2$ if $x\in A_3(r_3)\setminus A_2(r_2)$.  We can expand $g$ in a series
\[g(x)=g(\pi x +\rho x)=\sum _{m=0}^\infty g_m(\pi x + \rho x)\] where the $g_m$'s are $m$-homogeneous in $\rho x$, and convergence is uniform on compact subsets of $X$.  There exist $\delta > 0$ and $M\geq 0$ such that $\pi x\in A_3$ and $\| \rho x\| \leq \delta$ imply $\|g(x)\|\leq M$, hence
\begin{equation} 
	\|g_m(x)\|=\| \frac{1}{2\pi i}\int_{|\zeta|=\delta} \frac{g(\pi x +\zeta \rho x)}{\zeta^{m+1}} d\zeta\| \leq \delta^{-m}M. 
	\label{g_m}
\end{equation}

With $\alpha=\alpha_1 +\alpha_2$ and $\beta$ (sufficiently large) integers to be specified later, define functions from $A_4[r_4]$ to $Z_m=Y\overline{\otimes} W^{\overline{\otimes}(\alpha m +\beta)}$ by 
\[	h_m (x) =  g_m(x)\otimes \psi(x)^{\otimes(\alpha m +\beta)}.\]
Note that each $h_m$ is holomorphic.  We will show that $\tilde{h} = ( h_m)_{m=0}^\infty $ defines 	a holomorphic map $A_4(r_4)\to Z=\overline \bigoplus_{m=0}^\infty Z_m$.  Indeed, if $K\subset A_4(r_4)$ is compact and $x\in K$,   
\[ \|\tilde{h}(x)\| ^p \leq \sum_{m=0}^\infty (\|g_m(x)\|(\|\psi( x)\| ^{\alpha m +\beta})^p\]
\[\leq M \sum_{m=0}^\infty \|g_m(\pi x+M^\alpha \rho x)\|^p,\]
(here and below, we use $p=1$ in case (a)).  Now $\tilde{K}=\pi K + M^\alpha \rho K$ being compact and $g$ entire,  
\[ M \sum_{m=0}^\infty \|g_m(\pi x + M^\alpha \rho x)\|^p\] converges uniformly on $K$.  Thus $\tilde{h}$, being a limit of holomorphic functions which is uniform on compact sets, is holomorphic on $A_4(r_4)$. 

To estimate $\tilde h (x)$ for $x\in A_3(r_3)\setminus A_2(r_2)$ we use H\"{o}lder's inequality:
\begin{align*}
	 \| \tilde{h} (x)\| &=  (\sum _{m=0}^\infty (\| g _m(x)\| \| \psi( x)\| ^{\alpha m +\beta })^p)^{1/p}\\
	& \geq (\sum _{m=0}^\infty (\| g _m(x)\| 2^{\alpha m+\beta })^p )^{1/p}.\\
	& \geq (\sum _{m=0}^\infty \| g _m(x)\|)(\sum _{m=0}^\infty (1/2)^{(\alpha m +\beta )q})^{-1/q}
\end{align*}
  where $1/q+1/p=1$. If $\beta\geq 2$, then
	\[ \| \tilde{h}(x)\| \geq 2\sum _{m=0}^\infty  \| g _m(x)\| \geq 2\|g(x)\|\leq \|g(x)\|+\epsilon, \quad x\in A_3(r_3)/ A_2(r_2).\]

Next we estimate $\tilde{h}$ for $x\in A_1[r_1]$:  
\begin{align}
	\|\tilde{h}(x)\|&=(\sum _{m=0}^\infty (\| g _m(x)\|( \| \psi(x)\|)^{\alpha m +\beta })^p)^{1/p}\\
	&\leq ( \sum _{m=0}^\infty \| g_m(x) 2^{-(\alpha m +\beta )}\|^p)^{1/p}\\
	& \leq 2^{-\beta  } (\sum _{m=0}^\infty ( \| g_m(\pi x+2^{-\alpha_2} \rho x)\| 2^{-m\alpha_1 } )^p)^{1/p.} \label{h}
	\end{align} 
Choose $\alpha_2$ so that $2^{-\alpha_2}\|\rho x\| \leq \delta$ for $x\in A_1[r_1]$, then by (\ref{g_m}), we continue estimating (\ref{h}):
	\[  \leq 2^{-\beta  } (\sum _{m=0}^\infty (M \delta^{-m} 2 ^{-m\alpha_1}  )^p)^{1/p}. \]
	Now choose the exponent $\alpha_1$ so that $2^{-\alpha_1}\leq \delta/2$ and $\beta\geq2$ so that $2^{-\beta+1}M\leq \epsilon$, giving 
	\begin{align}  \|\tilde{h}(x)\|&\leq 2^{-\beta  }M (\sum _{m=0}^\infty (2 ^{-m-1})^p)^{1/p} \\
	& \leq \epsilon /2( \sum _{m=0}^\infty (2 ^{-m-1})^p)^{1/p} < \epsilon/2. \end{align}
Therefore $\|\tilde{h}\| < \epsilon /2$.
  Let $h\in \Oc(\Omega, Z)$ such that $h$ approximates $\tilde{h}$ within $\epsilon/2$ on $A_3(r_3)$.  Then $\|h(x)\|<\epsilon$ if $x\in A_1[r_1]$ and $\|h(x)\|> \|g(x)\|$ if $x\in A_3(r_3)\setminus A_2(r_2)$.
\end{proof}

\begin{prop} \label{dom}
Let $\epsilon>0$.  Assume either that (a) $X$ and $Y$ are Banach spaces with 1-unconditional basis or (b) $X$ and $Y$ are separable $L^p$ spaces.  Let $2^4\beta<\alpha<2^{-9}/c$ where $c$ is the unconditional basis constant of $X$.  If $\Omega\subset X$ is pseudoconvex, and $g\in \Oc(X;Y)$, then there are a Banach space $Z$ and $h\in \Oc (\Omega; Z)$ such that in case (a) $Z$ has 1-unconditional basis, in case (b) $Z$ is a separable $L^p$ space, and 
\begin{enumerate}
\item $\| h(x)\|_Z \leq \epsilon$ if $x\in \Omega _N \langle \beta \rangle$, and
\item $\|h(x)\|_Z \geq \|g(x)\|_Y$ if $x\in \Omega _{N+1} \langle \alpha \rangle \setminus \Omega_N\langle \alpha \rangle$.  
\end{enumerate}  
\end{prop}
 The proof of this theorem depends on showing that sets $A_i$ and functions $r_i$ can be defined satisfying the conditions of Lemma 6.2 so that $\Omega _N \langle \beta \rangle \subset A_1[r_1]$ and $\Omega _{N+1} \langle \alpha \rangle \setminus \Omega_N\langle \alpha \rangle \subset A_3(r_3)\setminus A_2(r_2)$.  Then we apply Lemma 6.2.  These containments are shown in the proof of Proposition 4.2 in [L2], here using $\mu=1/(2c)$.   

\begin{proof}[Proof of Proposition \ref{vb}]
The first part of the proof is as in [L2].  Let $\alpha<2^{-9}/c$, where $c$ is the unconditional basis constant of $X$.  First, for each $N\in \mathbb{N}$ we create a Banach space with a 1-unconditional basis (separable $L^p$ space) $Z_N$ and a function $g_N\in \Oc(X, Z_N)$ such that $u\leq \|g_N\|_{Z_N}$ on $\Omega_N\langle\alpha\rangle$.   Assume $u>1$ everywhere.  Define $A=\overline{\Omega_N\langle \alpha \rangle}\cap \pi_N X$, a compact set.  Recalling the definition of $d(x)$ from (\ref{d}), if $t\in A$, then 
\[ \Omega_N\langle \alpha \rangle \cap \pi_N^{-1}t \subset B(t,\alpha d(t)).\]
Therefore, there is a neighborhood $U\ni t$, $U\subset \pi_NX$, such that $\Omega_N\langle \alpha \rangle \cap \pi_N^{-1} U\subset B(t, 2\alpha d(t))$.  By compactness of $A$, there is a finite set $T\subset A$ such that 
\[ \Omega_N \langle \alpha \rangle \subset \bigcup_{t\in T} B(t, 2\alpha d(t)).\]
Now we show that $B_t =B(t, 4c\alpha d(t))\in \B$.  Indeed,  
\[ 2\text{ diam}B_t=16c\alpha d(t)<2d(t) \leq \text{diam}\Omega. \]
So by the hypothesis of Proposition \ref{vb}, there are a Banach space $V_t$ with 1-unconditional basis (separable $L^p$ space) and a function $f_t\in \Oc (B_t, V_t)$ such that $u\leq \|f_t\|_{V_t}$ on $B_t$.  By Lemma 4.3 concerning approximation on balls, there are $f'_t\in \Oc (X, V_t)$ such that $\|2f_t -f'_t\|_{V_t}< 1$ on $B(t, 2\alpha d(t))$; so $u\leq \|f'_t\|_{V_t}$ on $B(t, 2\alpha d(t))$.  Define $Z_N=\bigoplus _{t\in T} V_t$ with the $l^1$ norm in case (a) or $l^p$ norm in case (b), and $g_N=\bigoplus_{t\in T}f'_t \in \Oc(X, Z_N)$.  Then since the sets $B(t, 2\alpha d(t))$ cover $\Omega_N\langle \alpha \rangle$, given $x\in \Omega_N\langle \alpha \rangle$, there is a $t$ such that $\|f_t(x)\|\geq u(x)$, so $u(x)\leq \|g_N(x)\|_{Z_N}$ for all $x\in \Omega_N\langle \alpha \rangle$.

Fix $\beta <2^{-4}\alpha$.  By using Proposition \ref{dom}, we have $Y_N$ a Banach space with 1-unconditional basis (separable $L^p$ space) and $h_N\in \Oc (\Omega, Y_N)$ satisfying $\|h_N\|_{Y_N}\leq 2^{-N}$ on $\Omega_N\langle \beta \rangle$ and $\|h_N\|_{Y_N}\geq \|g_N\|_{Z_N}\geq u$ on $\Omega_{N+1}\langle \alpha \rangle \setminus \Omega_N\langle \alpha \rangle$.  Let $V=\overline\bigoplus_{N=0}^\infty  Y_N$.  Since any $x$ in $\Omega$ is in $\Omega_N\langle \beta \rangle$ for all but finitely many $N$, $\sum _N\|h_N\|$ converges locally uniformly, implying that  $f=\bigoplus _N h_N$ is in $\Oc(\Omega, V)$.  For any $x\in \Omega$, there is an $N$ so that $x\in \Omega_{N+1}\langle \alpha \rangle \setminus \Omega_N\langle \alpha \rangle$, then $\|f(x)\|_V\geq \|h_N(x)\|_{Y_N}\|\geq u(x)$,  and $\|f\|_V \geq u$ on $\Omega$.

\end{proof}

\end{document}